\documentclass[12pt]{article}

\usepackage[utf8]{inputenc}
\usepackage[T1]{fontenc}
\usepackage{lmodern}
\usepackage{microtype}
\usepackage{graphicx}
\usepackage{amsmath,amssymb,amsthm,mathrsfs,amsfonts,mathtools}
\usepackage{latexsym,xspace}
\usepackage{enumitem}
\usepackage[dvipsnames,svgnames]{xcolor}
\usepackage{todonotes}
\usepackage{comment}
\usepackage{tikz}
\usepackage{authblk}
\usepackage{booktabs}
\usepackage[a4paper,margin=1in]{geometry}
\usepackage{wrapfig}
\usepackage{parskip}
\usepackage{nicefrac}
\usepackage[font=small]{caption}
\usepackage{subcaption}
\usepackage[normalem]{ulem}
\usepackage[hidelinks]{hyperref}
\usepackage{cleveref}

\usepackage{tcolorbox}

\newtcolorbox{propertybox}{
  colback=blue!3!white,
  colframe=blue!45!black,
  boxrule=0.5pt,
  arc=2pt,
  left=8pt,
  right=8pt,
  top=6pt,
  bottom=6pt,
  before skip=8pt,
  after skip=8pt
}

\definecolor{requestred}{RGB}{185,28,28}

\usetikzlibrary{arrows.meta,fit,backgrounds}
\usetikzlibrary{positioning}
\usetikzlibrary{arrows}
\usetikzlibrary{decorations.pathmorphing}
\usetikzlibrary{decorations.markings}
\usetikzlibrary{decorations.pathreplacing}
\usetikzlibrary{shapes.geometric}
\tikzset{
    small circles/.style={circle,inner sep=2pt,fill=#1},
    hollow circles/.style n args={2}{circle,inner sep=#1,draw=#2,thick},
    stars/.style={star,inner sep=2pt}
}

\newcommand{\circlereq}[1]{%
  \tikz[baseline=(char.base)]\node[draw=requestred,text=requestred,circle,
  line width=.55pt,inner sep=.6pt,font=\bfseries\scriptsize] (char) {$#1$};}

\DeclareMathOperator{\E}{\mathbb{E}}
\newcommand{\arc}[2]{#1\mathrel{\scalebox{0.7}[1]{\(\to\)}}#2}

\theoremstyle{plain}

\newtheorem{theorem}{Theorem}
\newtheorem{lemma}[theorem]{Lemma}

\newtheorem{proposition}[theorem]{Proposition}

\newtheorem{claim}[theorem]{Claim}

\newcommand{\restatedtheoremname}{}
\newtheorem*{restatedtheoreminner}{\restatedtheoremname}

\newenvironment{restatedtheorem}[1]
  {%
    \renewcommand{\restatedtheoremname}{Theorem~\ref{#1}}%
    \begin{restatedtheoreminner}[restated]%
  }
  {%
    \end{restatedtheoreminner}%
  }

\theoremstyle{definition}

\theoremstyle{remark}

\newtheorem{property}{Property}

\newenvironment{claimproof}[1]{\par\noindent\textit{Proof of Claim:}\space#1}{\hfill $\blacksquare$}

\newcommand{\authoremail}[1]{%
  \href{mailto:#1}{\nolinkurl{#1}}%
}

\title{\textbf{List coloring degenerate graphs with requests}}

\author[1]{Rajat Adak}
\author[1]{L. Sunil Chandran}
\author[2]{Shiwali Gupta}
\author[3]{Tom\'{a}\v{s} Kaiser}
\author[2]{Subrahmanyam Kalyanasundaram}
\author[2]{Rogers Mathew}
\author[1]{Saraswati Girish Nanoti}
\author[1]{Suraj Kumar Sahoo}
\author[1]{Umesh Shankar}
\author[4]{Joel M. Thomas}

\affil[1]{Department of Computer Science and Automation, Indian
  Institute of Science, Bangalore, India}
\affil[2]{Department of Computer Science and Engineering, Indian
  Institute of Technology Hyderabad, India}
\affil[3]{Department of Mathematics and European Centre of Excellence
  NTIS (New Technologies for the Information Society), University of
  West Bohemia in Pilsen, Pilsen, Czech Republic}
\affil[4]{Department of Mathematics,  Indian Institute of Science
  Education and Research Tirupati, India}
\affil[ ]{Emails (in author order):
  \authoremail{rajatadak@iisc.ac.in}, \authoremail{sunil@iisc.ac.in},
  \authoremail{cs21resch11002@iith.ac.in}, 
  \authoremail{kaisert@kma.zcu.cz},  \authoremail{subruk@iith.ac.in}, 
  \authoremail{rogers@iith.ac.in}, 
  \authoremail{saraswatig@iisc.ac.in}, 
  \authoremail{surajks@iisc.ac.in}, 
  \authoremail{umeshshankar@outlook.com}, 
  \authoremail{joel_20221070@students.iisertirupati.ac.in}
}

\date{}

\date{November 2026}
\begin{document}
\maketitle
\begin{abstract}
  We prove three results on the flexible list coloring of
  $d$-degenerate graphs. Following Dvo\v{r}\'{a}k, Norin and Postle
  [J. Graph Theory 2019], a graph $G$ is $(k,\epsilon)$-flexible if
  for any assignment $L$ of lists of size $k$ to its vertices and any
  choice of requested colors by a set $\widehat V$ of the vertices
  from their respective lists, there is an $L$-coloring satisfying at
  least $\epsilon|\widehat V|$ of the requests. The complete graph on
  $d+1$ vertices shows that if all $d$-degenerate graphs are
  $(k,\epsilon)$-flexible for some list size $k$ and satisfaction
  ratio $\epsilon$, then $k \geq d+1$ and $\epsilon\leq\tfrac1{d+1}$.

  Our results concern values at or near these extremes. For list size
  $d+2$, we improve the satisfaction ratio from an inverse exponential
  to an inverse polynomial function of $d$. For satisfaction ratio
  $\tfrac1{d+1}$, we prove that it is achievable for list size $2d$ as well as
  $\left\lfloor(d+2)\left(2-\tfrac1e\right)\right\rfloor$ (which is
  less than $2d$ for $d\geq 9$).
\end{abstract}

\textbf{Brief History and Declaration of AI Assistance:} 
The third, fifth and sixth authors have been working on questions arising from the flexible list coloring problem for the last 3 years. Recently, in the summer of 2026, many of the authors worked on this problem as part of two visits at the Indian Institute of Science, Bangalore and Indian Institute of Technology, Hyderabad. 
The last author is on a year long visit to Indian Institute of Technology, Hyderabad as part of the work on his Masters thesis. 
Discussions with ChatGPT 5.6 Sol and 6.2 Sol were used. 
%More specifically, ChatGPT Sol 5.6 gave Theorem~\ref{thm:prime-power} when given \cite{DN19} as input. 
We had independently obtained an algorithmic version of \cite[Theorem 4]{KMMP24}. When prompted with this algorithm, ChatGPT provided  Theorem~\ref{thm:dcubed}. 
The counterexample in Proposition \ref{prop:chiflex-example}  was
discovered by ChatGPT.
% \rogers{
Theorem \ref{thm:e} was discovered by ChatGPT after persistent
prompting.
% }
This manuscript was written entirely by the authors, with AI tools
used only for minor editing tasks and drawing figures, and the authors take full
responsibility for the accuracy of its contents.

\section{Introduction}\label{intro}

A \emph{proper coloring} of a graph assigns colors to its vertices so that adjacent vertices receive distinct colors. The classical objective is to find such a coloring using as few colors as possible. The \emph{chromatic number} of a graph $G$, denoted $\chi(G)$, is the minimum number of colors required to do a proper coloring of $G$. Graph coloring has found applications in scheduling, frequency assignment in communication networks, register allocation in compilers, etc. \cite{chaitin1982register, lewis2021guide, marx2004graph}. To model situations in which the available resources vary across vertices, Vizing \cite{vizing1976vertex} and Erd\H{o}s, Rubin, and Taylor~\cite{erdos1979choosability} independently introduced \emph{list coloring}, also known as \emph{choosability}, in the 1970s. In this setting, each vertex is assigned a list of admissible colors, and the goal is to find a proper coloring in which every vertex receives a color from its list.

Although list coloring captures restrictions on admissible choices, it does not explicitly account for preferences among those choices. In applications such as scheduling and network design, certain admissible choices may be preferred because of cost, quality-of-service requirements, or user preferences. This motivates the study of \emph{flexible list coloring}, introduced by Dvo\v{r}\'{a}k, Norin and Postle \cite{DN19}. Here, each vertex in a specified subset has a preferred color from its list, and the objective is to find a proper list coloring that assigns the preferred color to as many of these vertices as possible. This can be viewed as a relaxation of the classical precoloring extension problem \cite{precoloringPartI-IntervalGraphs}.

In this paper, we obtain three results on flexible list coloring of degenerate graphs. Before we state them, let us recall the necessary definitions in detail. All graphs will be simple and undirected, and we generally follow the terminology of Diestel~\cite{Diestel}.

A \emph{list assignment} $L$ for a graph $G$ assigns a set $L(v)$ of admissible colors to each vertex $v \in V(G)$. An \emph{$L$-coloring} of $G$ is a function $f$ defined on $V(G)$ such that $f(v) \in L(v)$ for every vertex $v$. An $L$-coloring is \emph{proper} if adjacent vertices receive distinct colors, and $G$ is \emph{$L$-colorable} if it admits a proper $L$-coloring. A list assignment $L$ is a \emph{$k$-list assignment} if $|L(v)| = k$ for every $v \in V(G)$. The graph $G$ is \emph{$k$-choosable}, or \emph{$k$-list-colorable}, if it is $L$-colorable for every $k$-list assignment $L$. The smallest such integer $k$ is the \emph{list chromatic number} of $G$, denoted by $\chi_{\ell}(G)$. Taking identical lists at all vertices shows that $\chi(G) \leq \chi_{\ell}(G)$.

Given a list assignment $L$ for a graph $G$, a \emph{request function} for $G$ under $L$ is a function $r$ with a nonempty domain $\widehat{V} \subseteq V(G)$ such that $r(v) \in L(v)$ for every $v \in \widehat{V}$. For an $\epsilon \in (0,1]$, we say $G$ is \emph{$(k, \epsilon)$-flexible} if for every $k$-assignment $L$ and for every request function $r$ for $G$ under $L$ with domain $\widehat{V} \subseteq V(G)$, there exists a proper $L$-coloring $f$ of $G$ such that at least $\epsilon |\widehat{V}|$ vertices of $\widehat{V}$ satisfy $f(v) = r(v)$.
In this context, we refer to $k$ as the
\emph{list size} and to $\epsilon$ as the \emph{satisfaction ratio}.

This notion suggests two natural questions:
\begin{itemize}
\item For a given list size $k$, what is the largest $\epsilon$ for which $G$ is $(k,\epsilon)$-flexible?
\item For a prescribed satisfaction ratio $\epsilon$, what is the smallest $k$ that guarantees $(k,\epsilon)$-flexibility?
\end{itemize}

%\textcolor{red}{Comment by Subruk: The phrase ``satisfaction ratio'' is used at multiple places with a common meaning. Do we want to explicitly define it? Or leave it implicit as it is now? I prefer it to be explicitly defined the first time we use it in the introduction.}

Note that if $G$ is $(k,\epsilon)$-flexible for some $k$ and $\epsilon$, then $k \geq \chi_{\ell}(G)$. 

As mentioned above, we focus in this paper on flexible list coloring of degenerate graphs. Recall that a graph $G$ is \emph{$d$-degenerate} if each of its nonempty subgraphs has a vertex of degree at most $d$. It is easy to see that $G$ is $d$-degenerate if and only if it has a \emph{$d$-degeneracy ordering} of its vertex set, that is, an ordering $v_1,\dots,v_n$ such that each vertex $v_i$ has at most $d$ neighbours in $\{v_1,\dots,v_{i-1}\}$.
% A graph is $d$-degenerate if and only if there is an ordering of the vertices $v_1, v_2, \ldots, v_n$ such that every vertex $v_i$ has at most $d$ neighbors among $\{v_1, \ldots, v_{i-1} \}$. 
% The \emph{degeneracy} of $G$ is the smallest $d$ for which $G$
% is $d$-degenerate.
% We use $n$ for the number of vertices of a graph.
% We use the shorthand notation $uv$ for an edge $\{u,v\}$.

Whenever we prove all $d$-degenerate graphs to be $(k,\epsilon)$-flexible for some $k$ and $\epsilon$, we must have $k \geq d+1$ and $\epsilon \leq \tfrac1{d+1}$. Indeed, $(k,\epsilon)$-flexibility implies $k$-choosability, and there are $d$-degenerate graphs (such as $k_{d+1}$) of list chromatic number $d+1$. Furthermore, if all vertices of $K_{d+1}$ have the same list of size $k$ and request the same color, then no more than one of the requests can be granted, showing that $\epsilon \leq \tfrac1{d+1}$.

Dvo\v{r}\'{a}k, Norin and Postle \cite{DN19} showed that $d$-degenerate graphs are $(d+2,\frac{1}{(d+2)^{(d+1)^2}})$-flexible. The satisfaction ratio was improved by Kaul et al. \cite{KMMP24} who proved that $d$-degenerate graphs are $(d+2, \frac{1}{2^{d+1}})$-flexible. Their result was extended to weighted flexibility by Cambie and Cames van Batenburg~\cite{CC26}.

Our first result improves the denominator in the satisfaction ratio from exponential to polynomial in $d$. 
\begin{theorem}\label{thm:dcubed}
For every integer $d\ge0$, every $d$-degenerate graph is
\[
  \left(d+2,\frac{1}{(d+1)(d^2+d+1)}\right)\text{-flexible}.
\]
\end{theorem}

Kaul et al.~\cite{KMMP24} studied the relation between flexible list
coloring and `list packing'. Their Proposition 14 (together with the
fact that the `list packing number' of a $d$-degenerate graph is at
most $2d$, proved in~\cite{cambie2024packing}) implies that every
$d$-degenerate graph is $(2d, \frac{1}{2d})$-flexible. For this list size, we improve the satisfaction ratio to $\frac1{d+1}$. As noted above, this is the maximum value of $\epsilon$ we can hope for.

\begin{theorem}
\label{thm:2d-onebyd+1}
For any integer $d\ge 1$, every $d$-degenerate graph is
$$\left({2d, \frac{1}{d + 1}}\right)\text{-flexible}.$$
\end{theorem}

It turns out that for this extremal satisfaction ratio $\tfrac1{d+1}$, list size can be lowered when $d \geq 9$, as in our final main result.

\begin{theorem}\label{thm:e}
  For $d\geq 1$, every $d$-degenerate graph
  is
  \begin{equation*}
  \left(\left\lfloor (d+2)\left(2-\frac1e\right)\right\rfloor,\frac1{d+1}\right)\text{-flexible}.
  \end{equation*}
\end{theorem}

While the bound in Theorem~\ref{thm:e} is better than that in Theorem~\ref{thm:2d-onebyd+1} for any $d\geq 9$, we decided to include the latter as well, since its proof is shorter and simpler.

The above results say nothing about the other extremal situation for a $d$-degenerate graph $G$, namely when the list size is $d+1$. It is known \cite{DN19} that if $d+1$ is prime and $G$ has $n$ vertices, then $G$ is $(d+1, \tfrac{1}{n})$-flexible. Furthermore, for bipartite graphs $G$ this holds for every positive integer $d$~\cite{KMMP24}. On the other hand, we give in Section~\ref{sec:example} an example showing (for $d=2$) that a $d$-degenerate graph need not be $(d+1,\tfrac1{d+1})$-flexible.

Theorems~\ref{thm:dcubed} and \ref{thm:e} suggest the following open questions:
\begin{itemize}
\item Can the satisfaction ratio in Theorem~\ref{thm:dcubed} be further increased, say to $\Omega(\frac1{d^2})$?
\item Can the list size in Theorem~\ref{thm:e} be further decreased?
\end{itemize}

To conclude this section, we mention that many on the existing results on
flexible list coloring have focused on planar graphs;
see~\cite{dvovrak_triangle_free_planar_flexibility,
  dvorak_planar_girth_flexibility,yang_sufficient_conditions_planar_flexible}. Lidick\'{y}
et al.~\cite{LMMZ22} developed tools for weak flexibility, where every
vertex makes a request.

% the section below is by Rajat
% \section{$\left(d+2, \frac{1}{O(d^3)}\right)$-flexibility of $d$-degenerate graphs}
\section{Coloring from lists of size $d+2$}

In this section, we prove an $O(\tfrac1{d^3})$ upper bound for the satisfaction
ratio in a $d$-degenerate graph with a $d+2$-assignment.

%\note{we should put these in a proper environment}
%\rogers{RESPONSE: Done. Please check if this is fine.}
%\textcolor{red}{Comment by Subruk: I personally prefer the standard Theorem environment than this modified one. I can figure out a latex hack to preserve the numbering if that is OK with you, Rogers.}
%\textbf{Statement of Theorem \ref{thm:dcubed}:}
\begin{restatedtheorem}{thm:dcubed}
{For every integer $d\ge0$, every $d$-degenerate graph $G$ is
\[
  \left(d+2,\frac{1}{(d+1)(d^2+d+1)}\right)\text{-flexible}.
\]
}
\end{restatedtheorem}

\begin{proof} Fix a $d$-degeneracy ordering $v_1,\ldots,v_n$ of the
  vertices of $G$. Let $L$ be a $(d+2)$-list assignment and $r$ a
  request function with domain $\widehat{V} \subseteq V$.

For $i \in [n]$, let $N^-(v_i) = \{ v_j \in N(v_i) \mid j \in [i-1]\}$. For each $v\in \widehat{V}$, define
\begin{equation}\label{Udefn}
  U_v=\{r(v)\}\cup \bigcup_{\substack{x\in N^-(v)\\ r(v)\in L(x)}}L(x).
\end{equation}
Thus $U_v$ contains the requested color of $v$, that is $r(v)$, and the full list of
every neighbor of $v$ in $N^-(v)$ whose list contains $r(v)$. Since $|N^-(v)|\leq d$ and the list of each neighbor contains $d-1$ colors other than $r(v)$, we have
\begin{equation*}
  |U_v| \leq d^2+d+1.
\end{equation*}

Let $\mathcal P=\bigcup_{v\in V(G)}L(v)$ be the global palette of colors.
Choose a uniformly random total order $\prec$ of $\mathcal P$.  Call a vertex $v\in \widehat{V}$ \emph{good} when
\begin{equation*}
  r(v)=\min_{\prec}U_v,
\end{equation*}
that is, according to the total order, $r(v)$ is the first color to appear among the colors in $U_v$. 
Every element of $U_v$ is equally likely to be its minimum, and so
\begin{equation}\label{eq2}
      \Pr(v\text{ is good})=\frac1{|U_v|}\geq \frac{1}{d^2+d+1}.
\end{equation}

If $S$ is the set of good vertices, linearity of expectation and
\Cref{eq2} give
\begin{equation*}
  \mathbb E[S]
  =\sum_{v\in \widehat{V}}\frac1{|U_v|}
  \ge \frac{|\widehat{V}|}{d^2+d+1}.
\end{equation*}
Therefore some palette order satisfies
\begin{equation}\label{expectation}
  |S|\geq \frac{|\widehat{V}|}{(d^2+d+1)}.
\end{equation}
Note that $G[S]$ is also a $d$-degenerate graph, therefore $\chi(G[S]) \leq d+1$. Thus, there exists an independent set in $G[S]$, say $I$, such that \begin{equation}\label{mainbound}
    |I| \geq \dfrac{|\widehat{V}|}{(d+1)(d^2+d+1)}.
\end{equation}

Precolor every $v\in I$ with $r(v)$.  Since $I$ is an independent set, this precoloring is proper.

\begin{lemma}\label{lem:extend}
The precoloring $v\mapsto r(v)$ of $I$ extends to a proper $L$-coloring of $G$.
\end{lemma}

\begin{proof}
Start by coloring the uncolored vertices in the order
$v_1,\ldots,v_n$. When a vertex $v_i\notin I$ is reached, it has at
most $d$ neighbors in $N^-(v_i)$, all of which have already been
colored. Let
$$
L^-(v_i)
:=
\left\{
c\in L(v_i):
c\text{ is assigned to some vertex in }N^-(v_i)
\right\}.
$$
Since $|N^-(v_i)| \leq d$, we have
$|L^-(v_i)|\leq d.$

Now consider the precolored neighbors of $v_i$ that do not belong to
$N^-(v_i)$. Only those whose assigned colors belong to $L(v_i)$ can
reduce the number of colors available at $v_i$. If all such neighbors
are assigned the same color, then they collectively forbid at most one
additional color. Thus,
$$|L(v_i)\setminus L^-(v_i)|-1
\geq (d+2)-d-1
=1,$$
so at least one color remains available for $v_i$, and $v_i$ can be
colored greedily.

It therefore remains to rule out the possibility that there exist two
precolored vertices $x,y\in N(v_i)\setminus N^-(v_i)$,
such that $$r(x),r(y)\in L(v_i)
\qquad\text{and}\qquad
r(x)\neq r(y).$$

Since $x$ and $y$ are precolored, both belong to $I$. Since $x = v_j$ and $y = v_k$ where $j,k > i$, \Cref{Udefn} gives
\[
L(v_i)\subseteq U_x
\qquad\text{and}\qquad
L(v_i)\subseteq U_y.
\]

Since $x$ is good, $r(x)$ is the first color of $U_x$ in the fixed
total ordering. As $r(x)\in L(v_i)\subseteq U_x,$
it follows that $r(x)$ is also the first color of $L(v_i)$. Similarly,
since $y$ is good and $r(y)\in L(v_i)\subseteq U_y,$
the color $r(y)$ is the first color of $L(v_i)$ in the same total
ordering. A finite totally ordered set has a unique first element, and
therefore $r(x)=r(y)$, contradicting our assumption.

Hence the precolored neighbors of $v_i$ outside $N^-(v_i)$ forbid at
most one color of $L(v_i)$ in addition to the colors in $L^-(v_i)$.
Thus at most $d+1$ colors of the $(d+2)$-element list $L(v_i)$ are
forbidden, so an available color always remains. Coloring each
uncolored vertex in this manner produces a proper $L$-coloring
extending the precoloring of $I$.
\end{proof}

Thus, choose a palette order satisfying \Cref{expectation}.  By \Cref{lem:extend}, all requests on $I$, a set of order $\frac{|\widehat{V}|}{(d+1)(d^2+d+1)}$, can be satisfied
simultaneously in a proper $L$-coloring, which proves the theorem.
\end{proof}

\section{Satisfying $1/(d+1)$ of the requests with list size $2d$}
\label{sec:1-by-d+1-A}

As mentioned in Section~\ref{intro}, it follows from \cite{KMMP24} and
\cite{cambie2024packing} that every $d$-degenerate graph is
$(2d,\frac{1}{2d})$-flexible. We give an improvement of this result.

\begin{restatedtheorem}{thm:2d-onebyd+1}
{For any integer $d\ge 1$, every $d$-degenerate graph is 
\[ \left({2d, \frac{1}{d + 1}}\right)\text{-flexible}.\]}    
\end{restatedtheorem}

\begin{proof}
    Suppose $v_1,v_2, \ldots, v_n$ is a $d$-degeneracy ordering of the vertices of $G$, let $L$ be a $2d$-list assignment of $G$ and $r$ a request function. For $i\in [n]$, let $L(v_i)=\{w_{i,1}, w_{i,2}, \ldots ,w_{i,2d}\}$. Let $\widehat{V}$ be the domain of $r$.
    
    Suppose $\mathcal{C}$ is a set of $d+1$ proper $L$-colorings $C_1, C_2, \ldots , C_{d+1}$ of $G$ such that the following two properties always hold:
   
    \begin{propertybox}
    \begin{property}\label{prop:d+1-one}
        For any $v\in V(G)$ and for two distinct $1 \leq j,j'\leq d+1$, we have $C_j(v)\neq C_{j'}(v)$.
    \end{property}
    \begin{property}\label{prop:d+1-two}
        For all $v\in \widehat{V}$, there exists a coloring $C_j$ such that $C_j(v) = r(v)$, the requested color of $v$.
    \end{property}
    \end{propertybox}
    Then we call $\mathcal{C}$ a \emph{good $L$-coloring set of $G$}. We will first show that a good $L$-coloring set of $G$ always exists. 
    
    For $i\in [n]$, let $G_i$ be the subgraph of $G$ induced on the vertices $v_1,v_2,\ldots ,v_{i}$. We describe an iterative procedure that for every $i\in [n]$, produces a good $L$-coloring set $\mathcal{C}_i$ of $G_i$. As a consequence, $\mathcal{C}_n$ would be a good $L$-coloring set for $G$.

    A good $L$-coloring set $\mathcal{C}_1$ of $G_1$ clearly exists as we can color $v_1$ by arbitrarily picking $d+1$ colors in its list, always including its requested color in case $v_1\in \widehat{V}$. Now, suppose that we are on the $i$th iteration for some $i> 1$. By assumption, at the beginning of this iteration, we have already successfully obtained a good $L$-coloring set $\mathcal{C}_{i-1}$ of $G_{i-1}$. Let $\mathcal{C}_{i-1}=\{C_1, C_2, \ldots , C_{d+1}\}$. 
    We obtain a good $L$-coloring set $\mathcal{C}_{i}$ of $G_{i}$ by extending $\mathcal{C}_{i-1}$ in the following way.

    %Construct a bipartite graph $H_i$ with parts $X_i=\{x_{i,j }\colon j\in [d+1]\}$ and $Y_i=\{y_{i,k}\colon k\in [2d]\}$ such that $\{x_{i,j},y_{i,k}\}\in E(H_i)$ if and only if for every neighbor $u$ of $v_i$ in $G_i$, $C_j(u)\neq w_{i,k}$. Note that $C_j(u)$ is well defined for every neighbor $u$ of $v_i$ in $G_i$ because $u$ is a vertex of $G_{i-1}$. 

     %\textcolor{red}{Comment by Subruk: I am not sure if we needed sets $X_i, Y_i$ with vertices $x_{i,j}, y_{i,k}$. I have modified the writing  by using $\mathcal C_{i-1}$ and $L(v_i)$ instead of $X_i, Y_i$, which I found to be a notation that could be eliminated. }

    Construct a bipartite graph $H_i$ with parts $\mathcal C_{i-1}$ and $L(v_i)$ such that there is an edge between $C_j$ and $w_{i,k}$ if and only if for every neighbor $u$ of $v_i$ in $G_i$, $C_j(u)\neq w_{i,k}$. Note that $C_j(u)$ is well defined for every neighbor $u$ of $v_i$ in $G_i$ because $u$ is a vertex of $G_{i-1}$. 
    
    \begin{claim}\label{clm1}
       For all $j\in [d+1]$, the coloring $C_j$ has at least $d$ neighbors in $L(v_i)$.
    \end{claim}
    \begin{claimproof}
    Let $j\in [d+1]$. Let $S_{i,j}$ be the set of colors that appear in the neighborhood of $v_i$ in $G_i$ in the coloring $C_j$. The vertex $v_i$ has at most $d$ neighbors in the graph $G_i$ (because of the degeneracy ordering). 
    Hence, $|S_{i,j}|\leq d$. Therefore, $|L(v_i)\setminus S_{i,j}|\geq d$ and every color in $L(v_i)\setminus S_{i,j}$ is a neighbor of $C_j$.
    \end{claimproof}
    \begin{claim}\label{clm2}
        Every color in the list $L(v_i)$
        has at least 1 neighbor in $\mathcal C_{i-1}$. 
    \end{claim}
    \begin{claimproof}
        Suppose not. Let $k\in [2d]$ be such that $w_{i,k}$ has no neighbor in $\mathcal C_{i-1}$. Then for all colorings $C_j$, there is a neighbor $u$ of $v_i$ in $G_i$ such that $C_j(u)=w_{i,k}$.
        For each $j\in [d+1]$, choose a neighbor, say $u_j$, of $v_i$ in $G_i$ such that $C_j(u_j)=w_{i,k}$. Since $\mathcal{C}_{i-1}$ is a good $L$-coloring set of $G_{i-1}$, by Property~\ref{prop:d+1-one}, we must have that $u_1,u_2, \ldots , u_{d+1}$ are all pairwise distinct vertices. This contradicts the fact that $v_i$ has at most $d$ neighbors in $G_i$.
    \end{claimproof}
    \begin{claim}\label{clm3}
        The graph $H_i$ has a matching $M_i$ that saturates $\mathcal C_{i-1}$ such that if $v_i$ is a vertex of $\widehat{V}$, its requested color $r(v_i)$ is matched.
    \end{claim}
    \begin{claimproof}
        To show that $H_i$ has a matching that saturates $\mathcal C_{i-1}$, we need to check whether Hall's condition is satisfied, i.e., $\forall S\subseteq \mathcal C_{i-1}$, $|N_{\mathcal C_{i-1}}(S)|\geq |S|$. By Claim~\ref{clm1}, every coloring in $\mathcal C_{i-1}$ has degree at least $d$ in $L(v_i)$. Therefore, Hall's condition is satisfied
        for every subset of $\mathcal C_{i-1}$ of cardinality at most $d$.
        By Claim~\ref{clm2}, $N_{H_i}(\mathcal C_{i-1})= L(v_i)$. Since $|L(v_i)|=2d\geq d+1=|\mathcal C_{i-1}|$, Hall's condition is satisfied for the whole set $\mathcal C_{i-1}$ as well. Therefore, $H_i$ has a matching $M_i$ that saturates $\mathcal C_{i-1}$.

        Now suppose $v_i\in \widehat{V}$ and let $r(v_i)$ be its requested color. If $r(v_i)$ is a matched vertex in $M_i$, we are done. Suppose not. By Claim~\ref{clm2}, there is a neighbor $C_j\in \mathcal C_{i-1}$ of $r(v_i)$. Since, $M_i$ saturates $\mathcal C_{i-1}$, there is an edge incident on $C_j$ in $M_i$. Remove that edge from $M_i$ and add the edge $\{C_j,r(v_i)\}$ to $M_i$. After this substitution, $M_i$ is still a matching that saturates $\mathcal C_{i-1}$ and now has $r(v_i)$ as a matched vertex.  
    \end{claimproof}    
    
    Let $M_i$ 
    %$$M_i=\{\{x_{{i,1}},y_{i,k_1}\}, \{x_{{i,2}},y_{i,k_2}\}\ldots, \{x_{{i,d+1}},y_{i,k_{d+1}}\}\}$$ 
    be a matching as in the statement of Claim~\ref{clm3}, where each coloring $C_j\in \mathcal C_{i-1}$ is matched with a color $w_{i, k_j} \in L(v_i)$. We extend $\mathcal{C}_{i-1}$ to $\mathcal{C}_{i}$ by assigning $C_j(v_i)=w_{i,k_j}$ for every $j\in [d+1]$. By Claim~\ref{clm3}, $\mathcal{C}_{i}$ satisfies Properties \ref{prop:d+1-one} and \ref{prop:d+1-two}. Therefore, $\mathcal{C}_{i}$ is a good $L$-coloring set.

    Now pick a coloring $C_{j^*}\in \mathcal{C}_{n}$ that maximizes the number of vertices in $\widehat{V}$ that receive their requested color. Let $\widehat{V_1}, \widehat{V_2},\ldots , \widehat{V_{d+1}}\subseteq \widehat{V}$ be the sets of vertices that receive their requested colors in the colorings $C_1,C_2,\ldots , C_{d+1}\in \mathcal{C}_{n}$ respectively. Since $\mathcal{C}_{n}$ is a good $L$-coloring set of $G$, it is easy to verify that $\widehat{V}=\bigsqcup_{j\in [d+1]} \widehat{V_j}$. It follows that $|\widehat{V}_{j^*}|\geq \frac{1}{d+1}|\widehat{V}|$. Therefore, $G$ is $(2d,\frac{1}{d+1})$-flexible. 
\end{proof}

%%%%%%%%%%%%%%%%%%%%%%%%%%%%%%%%%%%%%%%%%%%%%%%%%%%%%%%%%%%%%%%%%%%%%%
% second part
%%%%%%%%%%%%%%%%%%%%%%%%%%%%%%%%%%%%%%%%%%%%%%%%%%%%%%%%%%%%%%%%%%%%%%

\section{Satisfying $1/(d+1)$ of the requests with smaller lists}
\label{sec:1-by-d+1-B}

For $d\geq 9$, the bound on list size in Theorem~\ref{thm:2d-onebyd+1} is not best possible, as shown by the following result.

\begin{restatedtheorem}{thm:e}
{For $d\geq 1$, every $d$-degenerate graph
  is
  \begin{equation*}
  \left(\left\lfloor (d+2)\left(2-\frac1e\right)\right\rfloor,\frac1{d+1}\right)\text{-flexible}.
\end{equation*}}
\end{restatedtheorem}

\begin{proof}
    Suppose that $v_1,v_2, \ldots, v_n$ is a $d$-degeneracy ordering of the vertices of $G$. For the integer $k = \left\lfloor (d+2)\left(2-\frac1e\right)\right \rfloor$, suppose $L$ is a $k$-list assignment of $G$. Let $r$ be a request function for $G$ with domain $\widehat{V}$. We extend the function $r$ to all of $V(G)$ by assigning each vertex outside $\widehat{V}$ an arbitrary color from its list.

    Our proof strategy is a generalization of the proof of Theorem \ref{thm:2d-onebyd+1}. Instead of maintaining $d+1$ colorings and extending them, we iteratively construct probability distributions on the valid colorings for the graph induced by 
    the first $i$ vertices. While Theorem~\ref{thm:2d-onebyd+1} ensures that one of the $d+1$ colorings satisfy a request, here we 
    ensure that the expected number of vertices in $\widehat{V}$ for which the requests are satisfied is at least $\frac{1}{d+1}|\widehat{V}|$.
    %Our proof strategy will be as follows. We will construct a probability distribution on all valid colorings and show that the expected number of vertices in $\widehat{V}$ for which the requests are satisfied is at least $\frac{1}{d+1}|\widehat{V}|$. 
    It follows that there is an assignment of colors that satisfies a $\frac{1}{d+1}$-fraction of the requests. The use of  probabilities instead of maintaining $d+1$ distinct colorings allows room for a better optimization of the list size needed.

    For $i\in [n]$, let $G_i$ be the subgraph of $G$ induced on the
    vertices $v_1,v_2,\ldots ,v_{i}$ and let $\mathcal C_i$ be the set of valid $L$-colorings of $G_i$. We describe an iterative
    procedure that creates the required probability distribution on
    $\mathcal C_i$ at step $i$. Set $M = \left\lfloor\frac{d+2}{e}\right\rfloor + 1$ and note that $k = 2d+4-M$. The reason
    for the choice of $M$ will be clear later. At each step $i$, we
    want our probability distribution to satisfy the following two
    properties:

\begin{propertybox}
    
\begin{property}\label{prop_1} For $c\in L(v_i)$, $\Pr[v_i \text{ is colored }c]\leq \frac{1}{d+1}$ with equality for $c=r(v_i)$.
\end{property}
\begin{property}
\label{prop_2}
For $c\in L(v_i)\setminus \{r(v_i)\}$ and a valid $L$-coloring $C$ of
$G_{i-1}$,
\begin{equation*}
  \Pr[v_i \text{ is colored }c\mid C]\leq \frac{1}{M}.
\end{equation*}
\end{property}
\end{propertybox}

To clarify the notation, the left hand side of the last inequality is the probability that a random $L$-coloring of $G_i$, sampled from the given distribution, assigns $c$ to $v_i$ conditioned on its agreement with $C$ on $G_{i-1}$.

For $i=1$, the required distribution is easy to find. The graph $G_1$ consists of a single vertex $v_1$. We assign the color $r(v_1)$ probability $\tfrac1{d+1}$, and divide the remaining probability among the other $2d+3-M$ colors in $L(v_1)$ equally. Thus, each of the other colors is assigned probability $\tfrac d{d+1}\cdot\tfrac1{2d+3-M}$. An elementary calculation shows that Properties~\ref{prop_1} and \ref{prop_2} hold. (Here, $G_0$ is the empty graph which has just one coloring.)

Let us therefore describe step $i$ of the construction, where $2 \leq i \leq n$, assuming that a probability distribution on $\mathcal C_{i-1}$ satisfying Properties \ref{prop_1} and \ref{prop_2} is given. The number $i$ will be fixed throughout the proof. We will obtain the required distribution on $\mathcal C_i$ by solving a flow problem in a network we define next.

    The flow network $H_i$ consists of the set $\mathcal C_{i-1}$ of all the possible colorings of $G_{i-1}$ that occur with non-zero probability, the set of colors $L(v_i)$, a source $s_i$ and a sink $t_i$. For $C\in \mathcal C_{i-1}$, define
    \begin{align*}
      A_C &= L(v_i)\setminus \{c:\,\text{color $c$ appears in }V(G_{i-1})\cap N(v_i)\text{ in the coloring }C\},\\
      \mu_C &= \Pr[C],
    \end{align*}
    where $\Pr[C]$ is the probability of the coloring $C$ in the distribution on $\mathcal C_{i-1}$.

     The source $s_i$ has an outgoing edge to every $C \in \mathcal C_{i-1}$. Every $C \in \mathcal C_{i-1}$ has outgoing edges to every color $c \in A_C$. Every $c \in L(v_i)$ has an outgoing edge to the sink $t_i$. These are the only edges of the flow network. The capacities of these edges are given as follows. The flow network is represented in Figure~\ref{fig:flow-network-cut}.
  $$\begin{array}{ll|c}
    \text{Edge}     &       & \text{Capacity} \\
    \hline
    \arc{s_i} C & (C \in \mathcal C_{i-1})     & \mu_C  \\
    \arc C c & (C\in \mathcal C_{i-1}, c\in A_C)  & \mu_C/M \\
    \arc c{t_i} & (c\in L(v_i))        & 1/({d+1})
  \end{array}$$

The sum of $\mu_C$ over all $C\in\mathcal C_{i-1}$ equals $1$. Indeed, it is the sum of all the probabilities in our probability distribution on $\mathcal C_{i-1}$. Consequently, all flows in the network $H_i$ have total value at most $1$, and each flow with total value $1$ assigns each edge $\arc{s_i} C$ value $\mu_C$. Our goal is to prove that indeed, a maximum flow in $H_i$ has total value $1$, and the following claim shows that it will provide us with the distribution we need. 

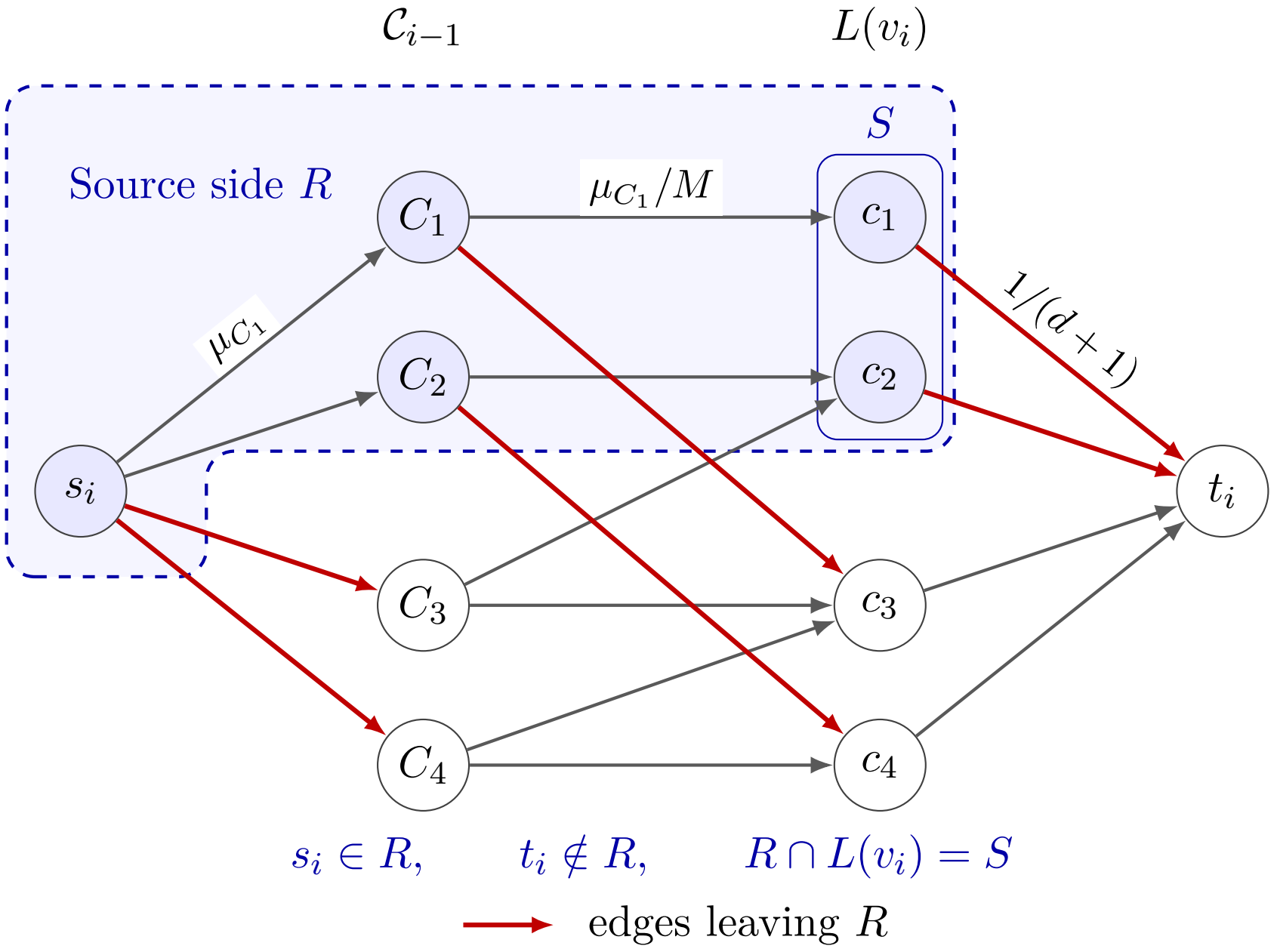
\begin{figure}[htbp]
\centering
\begin{tikzpicture}[
  x=1cm,y=1cm,font=\small,
  vertex/.style={
    circle,draw=black!75,fill=white,
    minimum size=8mm,inner sep=1pt
  },
  sourcevertex/.style={
    vertex,fill=blue!9
  },
  edge/.style={
    -{Latex[length=2mm]},draw=black!65,thick
  },
  cutedge/.style={
    edge,draw=red!75!black,very thick
  },
  capacity/.style={
    fill=white,inner sep=2pt,font=\footnotesize
  }
]

%------------------------------------------------------------
% Vertices in four layers.
% The upper two class vertices and color vertices lie in U.
%------------------------------------------------------------
\node[sourcevertex] (s) at (0,0)  {$s_i$};
\node[vertex]       (t) at (10,0) {$t_i$};

\foreach \j/\y in {1/2.4,2/1.0} {
  \node[sourcevertex] (C\j) at (3,\y) {$C_{\j}$};
  \node[sourcevertex] (c\j) at (7,\y) {$c_{\j}$};
}
\foreach \j/\y in {3/-1.0,4/-2.4} {
  \node[vertex] (C\j) at (3,\y) {$C_{\j}$};
  \node[vertex] (c\j) at (7,\y) {$c_{\j}$};
}

%------------------------------------------------------------
% Source side U: shaded region containing s_i, C_1, C_2,
% and precisely the color vertices in S.
%------------------------------------------------------------
\begin{scope}[on background layer]
  \path[
    draw=blue!65!black,fill=blue!4,dashed,
    rounded corners=7pt,thick
  ]
    (-0.65,-0.75) -- (1.1,-0.75) -- (1.1,0.35)
    -- (7.65,0.35) -- (7.65,3.55) -- (-0.65,3.55)
    -- cycle;
\end{scope}

\node[text=blue!65!black,anchor=west] at (-0.25,2.7)
  {Source side $R$};

% Layer labels.
\node at (3,4.05) {$\mathcal C_{i-1}$};
\node at (7,4.05) {$L(v_i)$};

% Mark S within the color layer.
\node[
  draw=blue!65!black,
  rounded corners=5pt,
  fit=(c1)(c2),
  inner sep=4pt,
  label={[text=blue!65!black]above:$S$}
] {};

%------------------------------------------------------------
% Source-to-class edges: capacity mu_C.
%------------------------------------------------------------
\draw[edge] (s) --
  node[capacity,above,sloped] {$\mu_{C_1}$} (C1);
\draw[edge] (s) -- (C2);

\foreach \j in {3,4}
  \draw[cutedge] (s) -- (C\j);

%------------------------------------------------------------
% Class-to-color edges: capacity mu_C/M.
% An edge C -> c exists exactly when c belongs to A_C.
%
% Illustrative choice:
% A_{C_1} = {c_1,c_3},   A_{C_2} = {c_2,c_4},
% A_{C_3} = {c_2,c_3},   A_{C_4} = {c_3,c_4}.
%------------------------------------------------------------
\draw[edge] (C1) --
  node[capacity,above] {$\mu_{C_1}/M$} (c1);

% Edges that do not leave U.
\foreach \a/\b in {2/2,3/2,3/3,4/3,4/4}
  \draw[edge] (C\a) -- (c\b);

% Edges from a class vertex in U to a color outside S.
\foreach \a/\b in {1/3,2/4}
  \draw[cutedge] (C\a) -- (c\b);

%------------------------------------------------------------
% Color-to-sink edges: capacity 1/(d+1).
% Exactly the edges from colors in S leave U.
%------------------------------------------------------------
\draw[cutedge] (c1) --
  node[capacity,above,sloped] {$1/(d+1)$} (t);
\draw[cutedge] (c2) -- (t);

\foreach \j in {3,4}
  \draw[edge] (c\j) -- (t);

%------------------------------------------------------------
% Cut condition and legend.
%------------------------------------------------------------
\node[text=blue!65!black] at (5,-3.2)
  {$s_i\in R,\qquad t_i\notin R,\qquad R\cap L(v_i)=S$};

\draw[cutedge] (3.35,-3.8) -- (4.15,-3.8);
\node[anchor=west] at (4.3,-3.8) {edges leaving $R$};

\end{tikzpicture}
\caption{The flow network $H_i$ and an $s_i$--$t_i$ cut with source
side $R$, where $R\cap L(v_i)=S$.
Red edges contribute to the cut capacity.
Capacities are shown on representative edges.}
\label{fig:flow-network-cut}
\end{figure}

\begin{claim}\label{claim:feasible}
    If a maximum flow in $H_i$ has total value $1$, then there is a probability distribution on $\mathcal C_i$ that satisfies Properties \ref{prop_1} and \ref{prop_2}.
\end{claim}
\begin{claimproof}
  Suppose that a flow $f$ in $H_i$ has total value $1$. As observed above, this implies $f(\arc{s_i}C) = \mu_C$, so $f$ must satisfy the following equality:
  \begin{equation}
    \label{eq:prelim-dist}
    \sum_{\substack{c'\in A_C}}f(\arc C {c'}) =\mu_C.
  \end{equation}
  In addition, just like any flow in $H_i$, $f$ clearly satisfies the following for each $C\in\mathcal C_{i-1}$ and $c\in L(v_i)$.
  \begin{eqnarray}
    f(\arc C c)&\le\frac{\mu_C}{M}, \label{eq:prelim-mu}\\
    \sum_{C'\in \mathcal C_{i-1}}f( \arc{C'}c)&\le \frac{1}{d+1}.\label{eq:prelim-marginal}
  \end{eqnarray}

  When defining the sought probability distribution on $\mathcal C_i$, we will use the following notation. If $C\in\mathcal C_{i-1}$ and $c\in A_C$, we write $C+c$ for the coloring that agrees with $C$ on $G_{i-1}$ and assigns $c$ to $v_i$. The symbol $\Pr[C+c]$ then denotes the probability assigned to $C+c$ by the new distribution on $\mathcal C_i$.

  It is a natural idea to define $\Pr[C+c]$ simply as $f(\arc C c)$. The inequality in Property~\ref{prop_1} would be satisfied by~\eqref{eq:prelim-marginal}. Furthermore, Property~\ref{prop_2} would hold by~\eqref{eq:prelim-mu} and the definition of conditional probability:
  \begin{equation*}
    \Pr[v_i \text{ is colored }c\mid C] = \frac{\Pr[C+c]}{\Pr[C]} = \frac{f(\arc C c)}{\mu_C} \leq \frac{\mu_C/M}{\mu_C} = \frac{1}{M}.
  \end{equation*}
Thus, the only possible problem is that $v_i$ might get its requested color $r(v_i)$ with probability less than $\tfrac1{d+1}$, in which case Property~\ref{prop_1} would fail. Define
  \begin{equation*}
    \beta \coloneqq  \sum_{C\in \mathcal C_{i-1}}f( C\rightarrow r(v_i)).
  \end{equation*}
  If $\beta = \tfrac1{d+1}$, then Property~\ref{prop_1} holds and we are done.

Suppose therefore that $\beta < \tfrac1{d+1}$ and define
\begin{align*}
  \alpha &\coloneqq \Pr[\text{no vertex in } V(G_{i-1})\cap N(v_i) \text{ is assigned color } r(v_i)],\\
  \theta&\coloneqq \frac{1}{\alpha-\beta}\left(\frac{1}{d+1}-\beta\right),
\end{align*}
where the probability in the first definition is with respect to the
distribution on $\mathcal C_{i-1}$.

Note that $0 < \theta \leq 1$. Indeed,
\begin{align*}
         \alpha &= 1-\Pr[r(v_i) \text{ appears in }V(G_{i-1})\cap N(v_i)]\\& \geq1- \sum_{u\in V(G_{i-1})\cap N(v_i)}\Pr[u\text{ is colored }r(v_i)]\geq1- \frac{d}{d+1}=\frac{1}{d+1}>\beta.
     \end{align*}
     Here, the first inequality follows by the union bound, and the second inequality follows by Property \ref{prop_1} at $u$.
     
     We now modify the definition of $\Pr[C+c]$ considered earlier. For every $C\in \mathcal C_{i-1}$ and $c\in A_C$, assign
     \begin{equation*}
       \Pr[C+c] \coloneqq
       \begin{cases}
         (1-\theta)\cdot f(\arc C c) + \theta \cdot \mu_C & \text{ if $c = r(v_i)$},\\
         (1-\theta)\cdot f(\arc C c) & \text{ if $c \neq r(v_i)$ and $r(v_i)\in A_C$},\\
         f(\arc C c) & \text{ if $r(v_i)\notin A_C$}.
       \end{cases}
     \end{equation*}
     This modification has the effect of boosting the probability of $v_i$ getting its requested color while maintaining the total probability. Since $\Pr[C+c]\leq f(\arc C c)$ whenever $c\notin r(v_i)$, Property~\ref{prop_2} and the inequality in Property~\ref{prop_1} remain valid for $c\notin r(v_i)$.

     It remains to show that Property~\ref{prop_1} holds at $v_i$ when $c = r(v_i)$.
     Note that $C\in \mathcal C_{i-1}$ has no edge in $H_i$ to $r(v_i)$ when $r(v_i)\notin A_C$. Therefore, we have
     \begin{align*}
       \Pr[v_i\text{ is colored }r(v_i)]
       &=
         \sum_{C\in \mathcal C_{i-1}}\Pr[ C+r(v_i)]\\
       &=\sum_{C\in \mathcal C_{i-1}\colon r(v_i)\in A_C} \Bigl((1-\theta) \cdot f(\arc C {r(v_i)})+\theta\cdot\mu_C\Bigr)\\
       &=
         (1-\theta)\cdot\beta+\theta \left(\sum_{C\in \mathcal C_{i-1}\colon r(v_i)\in A_C}\mu_C\right)\\
       &=
         (1-\theta)\cdot\beta+\theta \cdot \Pr[\text{no vertex in } V(G_{i-1})\cap N(v_i) \text{ is colored } r(v_i)]\\&=\beta-\theta \beta+\theta \alpha=\beta+\theta (\alpha-\beta)=\frac{1}{d+1}.
     \end{align*}
     
    % Moreover, it is easy to check that $\forall C\in \mathcal C_{i-1}$, $\sum_{c\in A_C}\Pr[C\wedge \{v_i\text{ is colored }c\} ]=\mu_C$.
     Hence, the probability distribution constructed above satisfies Properties \ref{prop_1} and \ref{prop_2} at the vertex $v_i$.
\end{claimproof}

Returning to the proof of Theorem~\ref{thm:e}, it remains to prove
that a maximum flow in $H_i$ has total value $1$. By the max-flow
min-cut theorem, this holds if and only if every $s_i-t_i$-cut has capacity at
least 1.

Consider a minimum $s_i-t_i$-cut in $H_i$; let us say it consists of
all edges between $R$ and $\overline R = V(H_i)\setminus R$ such that $s_i\in R$ and
$t_i\in \overline R$. We denote the cut by $\partial R$ and refer to $R$ and
$\overline R$ as its \emph{source side} and \emph{sink side},
respectively. See Figure~\ref{fig:flow-network-cut} for a
representation.

Let us denote $R\cap L(v_i)$ by $S$. The contribution made by
$C\in \mathcal C_{i-1}$ to the capacity of the cut is $\mu_C$ if it is
on the sink side and $\frac{\mu_C}{M}\cdot|A_C\setminus S|$ if it is
on the source side. By the minimality of $\partial R$,
$C$ belongs to $R$ or $\overline R$
depending on which of these values is smaller.
     %\rogers{Shouldn't it be `It follows that a good flow $f$ exists if and only if $\forall S\subseteq L(v_i)$,' } 
It follows that $H_i$ admits a flow of value $1$ if and only if
\begin{equation}\label{eq:feasible_flow}
         \frac{|S|}{d+1}+\sum_{C\in \mathcal C_{i-1}}\mu_C\min\left(1,\frac{|A_C\setminus S|}{M}\right) \geq 1,
       \end{equation}
     in which case we know there must be equality in view of the cut consisting of all edges incident with $s_i$.
       
     Let the neighbors of $v_i$ in $G_{i-1}$ be $v_{t_1},\dots,v_{t_m}$, and write $u_j = v_{t_j}$ for $j=1,\dots,m$. Since we are dealing with a $d$-degeneracy ordering, we have $m\leq d$.
     
     Consider a random coloring of $G_{i-1}$ as per the  probability distribution assigned to $\mathcal C_{i-1}$.
     For $j\in [m]$, let $U_j$ be the set of colors in $L(v_i)\setminus S$ that have not been assigned to $\{u_1, u_2, \ldots ,u_j\}$ and let $X_j=|U_j|.$ Observe that for each $j$,
     \begin{equation}
       \label{eq:Xj}
       X_j \geq k - j - |S|,
     \end{equation}
     because in the set $L(v_i)\setminus S$ of size $k-|S|$, at most $j$ colors are blocked by $U_j$.
     
     % \sout{Using this random variable $X_j$,}~\cref{eq:feasible_flow} \sout{can equivalently be written as,}
     Since $\mathbb{E}[X_m] = \sum_{C\in\mathcal C_{i-1}}\mu_C \cdot |A_C\setminus S|$, we can rewrite the left hand side of~\eqref{eq:feasible_flow} as
     \begin{equation*}
       \frac{|S|}{d+1} + \mathbb{E}\left[\min\left(1,\frac{X_m}{M}\right)\right].
     \end{equation*}
     Hence, in order to satisfy inequality \eqref{eq:feasible_flow}, it is enough to satisfy the following:
     \begin{equation}\label{eq:feasible_flow_equiv}
          \frac{|S|}{d+1} \geq \mathbb{E}\left[1-\min\left(1,\frac{X_m}{M}\right)\right].
     \end{equation}
     Note that the right hand side of~\eqref{eq:feasible_flow_equiv} is at most 1, so \eqref{eq:feasible_flow_equiv} holds whenever $|S|\geq d+1$.

     Therefore, we may assume $|S|\leq d$. 
For $j\in [m]$, let $Y_j=\min (M,X_j)$.
\begin{claim}\label{claim:Y_j_bound}
         For $j\geq 2$, $ \mathbb{E}[Y_j]\ge \left(1-\frac{1}{M}\right)\mathbb{E}[Y_{j-1}]-\frac{1}{d+1}.$
     \end{claim}
     \begin{claimproof}
         If $X_{j-1}\leq M$, then $Y_{j-1}=X_{j-1}$, and hence

$$Y_{j-1}-Y_j
=
\begin{cases}
1, & \text{if the color of }u_j\text{ is present in } U_{j-1},\\
0, & \text{otherwise}.
\end{cases}$$
Recall that $u_j=v_{t_j}$. For a coloring $C\in \mathcal C_{t_j-1}$, we get 
%\textcolor{red}{Comment by Subruk: I replaced $P_{t_j}$ by $\mathcal C_{t_j - 1}$ since I had replaced $P_i$ by $\mathcal C_{i-1}$. It needs to be checked if the meaning is correct.}
\begin{align*}
    &\mathbb{E}[Y_{j-1}-Y_j\mid C]\\&=\Pr[\text{the color of }u_j\text{ is present in } U_{j-1}\mid C]\\&
\le \Pr[\text{the color of }u_j\text{ is present in } U_{j-1}\setminus \{r(u_j)\}\mid C]+\Pr[u_j\text{ is colored }r(u_j)\mid C]\\&\le 
\frac{X_{j-1}}{M}
+
\Pr[u_j\text{ is colored }r(u_j)\mid C]\\& = %\le
\frac{Y_{j-1}}{M}
+
\Pr[u_j\text{ is colored }r(u_j)\mid C],
\end{align*}
where the second inequality follows by an application of Property~\ref{prop_2} at $v_{t_j}$.
On the other hand, if $X_{j-1}>M$, then
$X_{j-1}\geq M+1$ and we have
$X_j\geq M$. Consequently, in this case $Y_j=Y_{j-1}=M$, and hence, $Y_j-Y_{j-1}=0$.

\noindent 
Combining these two cases, for $C\in \mathcal C_{t_j-1}$, we have $$\mathbb{E}[Y_{j-1}-Y_j\mid C]\le \frac{Y_{j-1}}{M}
+
\Pr[u_j\text{ is colored }r(u_j)\mid C].$$
Using the law of total expectation, $$\mathbb{E}[Y_{j-1}-Y_j]\le
\frac{\mathbb{E}[Y_{j-1}]}{M}
+
\Pr[u_j\text{ is colored }r(u_j)] = \frac{\mathbb{E}[Y_{j-1}]}{M}
+ \frac{1}{d+1},$$
where the equality follows by an application of Property~\ref{prop_1} at $u_j = v_{t_j}$. A rearrangement of the terms gives us the claimed statement.
     \end{claimproof}
  
We now prove inequality~\eqref{eq:feasible_flow_equiv}. We know from~\eqref{eq:Xj} that $X_m \geq k-m-|S|$. Therefore, if $m \leq k-|S|-M$, then $X_m \geq M$ and the right hand side of~\eqref{eq:feasible_flow_equiv} is zero,
proving the inequality.
% }

%   \textcolor{red}{
We may thus assume that $m\ge k-|S|-M+1$.
Set $\ell=d-k+|S|+M$. Since by~\eqref{eq:Xj}, $Y_{k-|S|-M}=M$, and since $m\le d$,
we may apply Claim~\ref{claim:Y_j_bound} at most $\ell$ times to obtain
a lower bound on $\E[Y_m]$: %This bound decreases with the number
%of iterations, so
\begin{align*}
  \E[Y_m]
  &\ge
  M\left(1-\frac{1}{M}\right)^\ell
  -\frac{1}{d+1}
   \sum_{j=0}^{\ell-1}\left(1-\frac{1}{M}\right)^j\\
  &=
  M\left(1-\frac{1}{M}\right)^\ell
  -\frac{M}{d+1}
   \left[1-\left(1-\frac{1}{M}\right)^\ell\right].
\end{align*}
Consequently,
\[
  \E\left[1-\min\left(1,\frac{X_m}{M}\right)\right]
  =1-\frac{\E[Y_m]}{M}
  \le
  \frac{d+2}{d+1}
  \left[1-\left(1-\frac{1}{M}\right)^\ell\right].
\]
To show inequality~\eqref{eq:feasible_flow_equiv},  it suffices to prove that
\begin{equation}\label{eq:analytic}
  (d+2)\left[1-\left(1-\frac{1}{M}\right)^\ell\right]
  \le |S|.
\end{equation}
% }

% \textcolor{red}{
%Since $|S|-\ell=k-d-M$, we seek an upper bound on the
%left-hand side of~\cref{eq:analytic} minus $\ell$ that is
%independent of $S$.
The inequality
$\log\left(1+\frac{1}{M-1}\right)\le\frac{1}{M-1}$
gives
\[
  1-\frac{1}{M}
  =\frac{M-1}{M}
  \ge e^{-1/(M-1)}.
\]
Therefore,
\[
  (d+2)\left[1-\left(1-\frac{1}{M}\right)^\ell\right]-\ell
  \le
  (d+2)\left(1-e^{-\ell/(M-1)}\right)-\ell.
\]
% }

% \textcolor{red}{
To bound the expression on the right, observe that for $a,x>0$,
the function
\[
  g(z)=a\left(1-e^{-z/x}\right)-z,
  \qquad z\in\mathbb R,
\]
attains its maximum at $z=x\ln(a/x)$, with value
\[
  a-x-x\ln(a/x).
\]
We apply this with $a=d+2$ and $x=M-1$.
% }
%Substituting $x=M-1$ and $a=d+2$, 

The attained maximum value is $d+3-M-(M-1)\ln(\frac{d+2}{M-1})$. Let this value be denoted by $\nu$. If we show that $\nu\le k-M-d$, then we are done because $\nu+\ell\le  k-M-d+\ell=|S| $. 

We have $\nu\le k-M-d$ if and only if $k\geq 2d+3-(M-1)\ln(\frac{d+2}{M-1})$. Substituting $z$ for $M$ and letting $z$ vary, the right hand side is smallest when $z=\frac{d+2}{e}+1$ because $x\ln( \frac{a}{x})$ attains its maximum at $x=a/e$. At this extreme, the right hand side equals $2d+3-\tfrac{d+2}e$.

Therefore, when $k=2d+3-\left\lfloor\frac{d+2}{e}\right\rfloor = \left\lfloor (d+2)\left(2-\frac1e\right)\right \rfloor $ and $M=\left\lfloor\frac{d+2}{e}\right\rfloor+1$, inequality~\eqref{eq:feasible_flow_equiv} holds and the maximum flow in the network $H_i$ has total value $1$. By~Claim~\ref{claim:feasible}, we get a probability distribution on $\mathcal C_i$ that satisfies Properties \ref{prop_1} and \ref{prop_2}. %\note{(label them more clearly)} \rogers{(RESPONSE: calling them properties now and they are highlighted)}. 

Performing the construction for $i=1,\dots,n$, we eventually obtain a distribution on $\mathcal C_n$ that satisfies Properties \ref{prop_1} and \ref{prop_2}. Let $t$ denote the number of vertices of $\widehat{V}$ that are satisfied in a random coloring of $G_n$ under the distribution constructed above. Then $$\E[t]=\sum_{v\in \widehat{V}}\Pr[v\text{ is colored }r(v)]=\frac{ |\widehat{V}|}{d+1}\;.$$ 
This ensures that there is a coloring that satisfies the requests of at least a $\frac{1}{d+1}$-fraction of the vertices of $\widehat{V}$. Hence, $G$ is $\left(\left\lfloor (d+2)\left(2-\frac1e\right)\right \rfloor, \frac{1}{d+1}\right)$-flexible.
\end{proof}

%%%%%%%%%%%%%%%%%%%%%%%%%%%%%%%%%%%%%%%%%%%%%%%%%%%%%%%%%%%%%%%%%%%%%%

%the following section is due to Suraj
\section{List size $d+1$ and satisfaction ratio $\tfrac1{d+1}$}
\label{sec:example}

In this section, we give an example (for $d=2$) of a $d$-degenerate graph that is not $(d+1, \frac{1}{d+1})$-flexible. The example implies an affirmative answer to \cite[Question 10]{KMMP24} concerning the `list flexibility number' of degenerate graphs.  

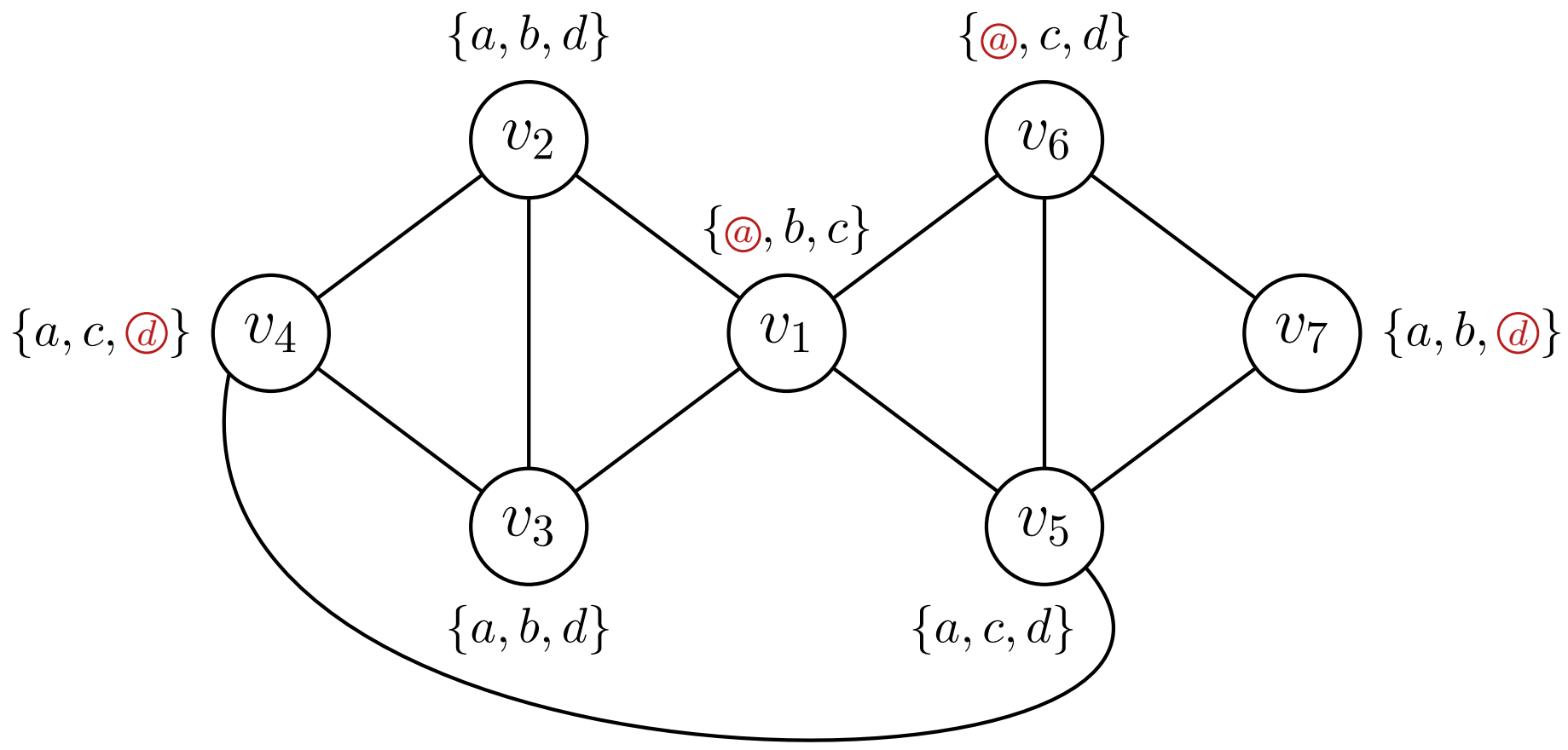
\begin{figure}[htbp]
\centering
\begin{tikzpicture}[
    vertex/.style={
        circle, draw, thick,
        minimum size=9mm,
        inner sep=0pt,
        font=\large
    },
    every label/.style={font=\small},
    line cap=round,
    line join=round
]
% Two congruent diamonds sharing v_1
\node[vertex, label=above:{$\{\circlereq{a},b,c\}$}]
    (v1) at (0,0) {$v_1$};

\node[vertex, label=above:{$\{a,b,d\}$}]
    (v2) at (-2,1.5) {$v_2$};

\node[vertex, label=below:{$\{a,b,d\}$}]
    (v3) at (-2,-1.5) {$v_3$};

\node[vertex, label=left:{$\{a,c,\circlereq{d}\}$}]
    (v4) at (-4,0) {$v_4$};
    
\node[vertex, label={[xshift=-4mm]below:{$\{a,c,d\}$}}]
    (v5) at (2,-1.5) {$v_5$};

\node[vertex, label=above:{$\{\circlereq{a},c,d\}$}]
    (v6) at (2,1.5) {$v_6$};

\node[vertex, label=right:{$\{a,b,\circlereq{d}\}$}]
    (v7) at (4,0) {$v_7$};

% Left diamond and its vertical diagonal
\draw[thick] (v1) -- (v2) -- (v4) -- (v3) -- (v1);
\draw[thick] (v2) -- (v3);

% Right diamond and its vertical diagonal
\draw[thick] (v1) -- (v6) -- (v7) -- (v5) -- (v1);
\draw[thick] (v6) -- (v5);

% Remaining edge, curved below the vertex labels
\draw[thick]
    (v4.south west)
    .. controls (-5,-3.8) and (4,-3.8) ..
    (v5.south east);

\end{tikzpicture}
 \caption{A 2-degenerate graph with flexible chromatic number at least 4. The list of colors are from the set $\{a,b,c,d\}$ and each vertex receives a list of size 3. The requested colors for the vertices $v_1,v_4,v_6$ and $v_7$ are circled in red. }
    \label{fig:chiflex-example}
\end{figure}

\begin{proposition}\label{prop:chiflex-example}
    There is a $2$-degenerate graph $G$ on $7$ vertices that is not $(3,\frac{1}{3})$-flexible.
%    and $\chi_{\mathrm{flex}}(G)\geq 4$. 
\end{proposition}
\begin{proof}
Consider the graph $G$ shown in Figure~\ref{fig:chiflex-example}. 
%\Subruk{Should we describe the graph in text? Or is the image enough? ANS: (Rogers). I think the figure is enough}
It is clear that in the ordering $v_1,v_2,v_3,v_4,v_5,v_6,v_7$, every vertex has at most 2 neighbors to its  left. Therefore, $G$ is 2-degenerate. The lists $L(v_i)$, for each $1 \leq i \leq 7$, are shown in Figure \ref{fig:chiflex-example}. Let $\widehat{V} := \{v_1, v_4, v_6, v_7 \}$ and $r~:~\widehat{V}\rightarrow \{a,d\}$ be defined as $r(v_1) = r(v_6) = a$ and $r(v_4) = r(v_7) = d$. Thus, $\left\lceil\frac{|\widehat{V}|}{3}\right\rceil =\left\lceil\frac{4}{3}\right\rceil=2 $. We show that under the given list assignment,  there is no proper list coloring of $G$ that satisfies at least 2 of the requests. This ensures that $G$ is not $(3, \frac{1}{3})$-flexible. 

Suppose $v_1$ and $v_4$ are given their requested colors, i.e., $a$ and $d$ respectively. Then the vertices $v_2$ and $v_3$ are only left with $b$ as their option and hence, cannot be colored. Similarly, if $v_1$ and $v_7$ are given their requested colors, then $v_5$ and $v_6$ are only left with $c$ as their option. And, $v_1$ and $v_6$ cannot be given their requested colors because they are adjacent and have the same requested color $a$. If $v_4$ and $v_6$ (or $v_7$) are given their requested colors, then $v_1$ is forced to $b$ and hence both $v_2$ and $v_3$ are only left with $a$ as their option. Finally, if $v_6$ and $v_7$ are given their requested colors, then $v_5$ and $v_1$ are forced to $c$ and $b$ respectively and this forces the triangle $v_2v_3v_4$ to be colored with 2 colors, namely $a$ and $d$, which is not possible.
\end{proof}

%%%%%%%%%%%%%%%%%%%%%%%%%%%%%%%%%%%%%%%%%%%%%%%%%%%%%%%%%%%%%%%%%%%%%%
% acknowledgments and refs
%%%%%%%%%%%%%%%%%%%%%%%%%%%%%%%%%%%%%%%%%%%%%%%%%%%%%%%%%%%%%%%%%%%%%%

\noindent \textbf{Acknowledgments:} The authors would like to thank Kunal Bag and Nitin Saurabh for helpful discussions. The fourth author thanks the Department of Computer Science and Automation at IISc Bangalore for its hospitality during his visit in July 2026. The fifth author thanks Louis Esperet for hosting him at Grenoble during the summer of 2025 and working on this problem. The fifth author likes to thank ANRF India for supporting this work through grant ANRF/ARGM/2025/002583/MTR. The seventh author likes to thank ANRF India for supporting this work through the NPDF fellowship research grant SP/ANRF 25-0599.06.

\bibliographystyle{plain}
\bibliography{bibfile}
\end{document}